\magnification\magstep1 \input amstex \input amsppt.sty


\font\bbf=cmbx10 scaled \magstep1   
\def\cal#1{{\Cal#1}}

\def\Pows{\Cal P\kern-.4mm_s\kern.3mm}
\def\lei{\hbox{\kern.45mm$_{^\downarrow}\kern-1.25mm\cap\kern.85mm$}}
\def\Qe{\hbox{$Q\kern-2.6mm\raise.2mm\hbox{\font\=cmssqi8\I}\kern1.7mm$}} 

\def\inve{\lower.85mm\hbox{$^{^-}$}\kern-.5mm{}^\iota}
\def\Alpha{\hbox{\font\=cmmi10 scaled\magstep1\\char'013}\kern0.15mm}
\def\Eps{\hbox{\font\=cmmi10 scaled\magstep1\\char'017}\kern0.15mm}
\def\Nu{\hbox{\font\=cmmi10 scaled\magstep1\\char'027}\kern0.15mm}
\def\Tau{\hbox{\font\=cmmi10 scaled\magstep1\\char'034}\kern0.15mm}
\def\value{\hbox{\kern.2mm\font\=cmr10\\char'022\kern-.2mm}} 
\def\image{\hbox{\font\=cmr10\\char'022\kern-1mm\char'022}} 
\def\images{\hbox{\font\=cmr10\\char'022\kern-1mm\char'022\kern-1mm\char'022}} 

\def\Circ{\kern.8mm\hbox{\font\=cmbsy10\\char'016}\kern.5mm}
\def\Examplee{{\font\=cmssi10\E\kern.15mmx\kern.15mma\kern.15mmm\kern.14mmp\kern.17mml\kern.15mme}\kern.3mm. }
\def\Examples{{\font\=cmssi10\E\kern.15mmx\kern.15mma\kern.15mmm\kern.14mmp\kern.17mml\kern.15mme\kern.15mms}\kern.3mm. }
\def\Remarkk{{\font\=cmssi10\R\kern.15mme\kern.15mmm\kern.15mma\kern.15mmr\kern.15mmk}\kern.3mm. }
\def\Remarkss{{\font\=cmssi10\R\kern.15mme\kern.15mmm\kern.15mma\kern.15mmr\kern.15mmk\kern.15mms}\kern.3mm. }
\def\N{{I\!\!N}} 
\def\No{{I\!\!N\kern-.54mm\lower.15mm\hbox{$_{\roman o}$}}} 
\def\Nopot#1{I\!\!N\kern-.54mm\lower.15mm\hbox{$_{\roman o}$}\kern-.7mm{}^{#1}} 
\def\potNo{^{\kern.37mm I\!\!{N_{}}_{\kern-.22mm\roman o}}} 
\def\Re{{I\!\!R}}

\def\plusinfty{\lower1.05mm\hbox{$^+$}\infty}
\def\minusinfty{\lower1.05mm\hbox{$^-$}\infty}
\def\dualbeta{^{\kern0.4mm\prime}_{\kern-.15mm\raise.95mm\hbox{$_{_\beta}$}}} 
\def\dualbetaa{^{\kern0.4mm\prime}_{\kern-.2mm\raise.95mm\hbox{$_{_\beta}$}}}
\def\duaL#1beta{^{\kern0.4mm\prime}_{\kern-.2mm\kern#1mm\raise.95mm\hbox{$_{_\beta}$}}}
\def\LL^#1{L\kern0.15mm\raise.4mm\hbox{$^{#1}$}\kern0.15mm}
\def\Ce{{\hbox{$C\kern-2.5mm\raise.2mm\hbox{\font\=cmssqi8\I}\kern1.48mm$}}}
\def\imag{\kern.15mm\lower.6mm\hbox{$^{^*}$}\kern-1.8mm\imath\kern.1mm} 

\def\biit#1{\hbox{\font\=cmmib10\#1}} 
\def\bmii#1#2{\hbox{\font\=cmmib#1\#2}}

\def\bcal#1#2{\hbox{\font\=cmbsy#1\#2\kern.3mm}}%
\def\Frak#1{\hbox{$\frak#1\kern.3mm$}}
\def\fssi#1{\hbox{\font\=cmssi10\#1}\kern0.15mm} 
\def\smb#1{\hbox{\font\†=cmmi8\†#1\kern.3mm}} 
\def\ssmb#1{\hbox{\font\=cmmi6\#1}} 
\def\ecal#1{\kern.1mm\hbox{\font\†=cmsy8\†#1\kern.3mm}} 
\def\ncal#1{\kern.1mm\hbox{\font\†=cmsy9\†#1\kern.3mm}} 
\def\vcal#1{\kern-.1mm\vec{\kern.2mm\hbox{\font\†=cmsy7\†#1}\kern.3mm}} 
\def\conc{\!\bold{\hat{\phantom w}}\!}

\def\idv{\hbox{\font\=cmr10\id}\kern.25mm\lower.8mm\hbox{\font\=cmr7\v}\kern.3mm} 
\def\idm{\hbox{\font\=cmr10\id}\kern.25mm\lower.8mm\hbox{\font\=cmr7\m}\kern.3mm} 
\def\seq#1{\langle#1\rangle}

\def\SemiNor{\Cal S_{_N}\kern0.15mm}
\def\vecs{\upsilon\kern-0.3mm\lower.15mm\hbox{$_s$}\kern0.2mm} 
\def\vecss{\hbox{\font\=cmitt10\v}\kern-0.1mm\lower.15mm\hbox{$_s$}\kern0.2mm} 
\def\bnull#1{\hbox{\font\=cmssbx10\0}{}_{\font\=cmmi6\lower.15mm\hbox{\kern-.1mm\#1\kern.15mm}}} 
\def\bzero#1{\hbox{\font\=cmbx10\0}{}_{\font\=cmmi6\lower.15mm\hbox{\kern-.1mm\#1\kern.15mm}}} 
\def\dom{{{}^{}\roman{dom}\,{}_{{}^{}}}}
\def\domr#1{\roman{dom}^{\font\=cmr6\raise.0mm\hbox{\kern.3mm\#1}}}

\def\rng{{}^{}\roman{rng}\,{}_{{}^{}}}
\def\CPi#1{C\kern-.2mm\lower.05mm\hbox{$_{_\Pi}$}\kern-1.52mm{}^{#1}}

\def\Cinfty{C\kern.4mm\raise.3mm\hbox{$^\infty$}\kern.15mm}
\def\Cinftyzero{\hbox{$C\kern.4mm\raise.3mm\hbox{$^\infty$}\kern.15mm\kern-3.5mm_{\font\=cmr6\lower.15mm\hbox{\kern.1mm\0}}\kern1.9mm$}}

\def\RHB#1#2{\raise#1mm\hbox{$#2$}} 
\def\LHB#1#2{\lower#1mm\hbox{$#2$}} 
\def\sssr^#1_#2{^{\,#1}_{\phantom\iota\kern-1.1mm\raise.5mm\hbox{$_{_{#2}}$}}} 
\def\ai#1{{}_{\font\=cmmi6\lower.15mm\hbox{\kern-.1mm\#1\kern.15mm}}} 
\def\yi#1{^{\font\=cmmi6\raise.0mm\hbox{\kern-.1mm\#1\kern.15mm}}} 
\def\ar#1{{}_{\font\=cmr6\lower.15mm\hbox{\kern.1mm\#1}}} 
\def\yr#1{^{\font\=cmr6\raise.0mm\hbox{\kern.3mm\#1}}} 
\def\lupar{\kern.2mm\lower1mm\hbox{$^{^(}$}} 
\def\rupar{\lower1mm\hbox{$^{^)}$}\kern-.15mm} 
\def\yyi#1{^{\font\=cmmi6\lower.6mm\hbox{\kern-.25mm\#1\kern-.05mm}}} 
\def\yyr#1{^{\font\=cmr6\lower.45mm\hbox{\kern-.25mm\#1\kern-.15mm}}} 
\def\yplus{\lower1mm\hbox{$^{^+}$}} 
\def\yminus{\lower1mm\hbox{$^{^-}$}} 
 %
\def\ClT{\roman{Cl}\kern.25mm\lower.4mm\hbox{$_{\Cal T}$}\kern0.2mm} 
\def\IntT{\sp\roman{Int}\kern.2mm\lower.4mm\hbox{$_{\Cal T}$}\kern0.2mm} 
\def\inc{\subseteq}

\def\exi#1{\exists\,#1\kern.2mm\,;}
\def\all#1{\forall\,#1\kern.2mm\,;}
\def\imply{\Rightarrow}

\def\sp{\kern0.15mm} 
\def\ssp{\kern0.37mm} 
\def\sssp{\kern0.52mm} 
\def\snn{\kern-0.2mm} 
\def\sn{\kern-0.3mm} 
\def\ssn{\kern-0.63mm} 
\def\KP#1{\kern#1mm} 
\def\KN#1{\kern-#1mm} 

\def\NS{\vskip1.7mm}
\def\NSN{\vskip1.7mm\noindent}
\def\œ$#1${\hbox{$#1$}} 
\def\"{\"a} \def\"{\"o}
\def\newProCla#1\par#2\par{\vskip1.7mm\noindent\bf#1\it#2\vskip1.7mm}
\def\Prooff{{\font\=cmssi10\P\kern.37mmr\kern.37mmo\kern.37mmo\kern.37mmf\kern.37mm. }\rm}

\def\QED{\hfill\hbox{$\ \sqcap\kern-2.45mm\sqcup$}}
\def\newQED{\hfill\hbox{$\ \sqcap$\hskip-2.45mm$\sqcup$}\vskip1.7mm}
\def\abstract#1{{\eightpoint\parindent5mm\narrower\baselineskip3.36mm#1\par}}
\def\noin{\noindent}
\def\Newline{\kern-10mm\newline}
\font\rp=cmr8

\def\RunMyHead#1#2#3#4{%
 \headline{\ifnum\pageno=\firstpage\hfil%
           \else{\ifodd\pageno{\rp#3\phantom\folio\hfil#4\hfil\phantom{#3}\folio}%
                 \else{\rp\folio\phantom{#2}\hfil#1\hfil\phantom\folio#2}%
                 \fi}%
           \fi}%
 \footline{\ifnum\pageno=\firstpage\hfil{\rp[\,\folio\,]}\hfil%
           \else\hfil%
           \fi}%
}%
\def\subhead#1\par#2\par{\vskip4mm\smallbreak\null\smallskip\vbox{\noindent\bbf#1\hfill\kern1.5mm#2\hfill\phantom{#1}\vskip2.5mm\nopagebreak}\nopagebreak\noindent}
\def\subheadd#1\par#2\par#3\par{\vskip4mm\smallbreak\null\smallskip\vbox{\noindent\bbf#1\hfill#2\hfill\phantom{#1}\vskip1.5mm\centerline{#3}\vskip2.5mm\nopagebreak}\nopagebreak\noindent}
\def\CinftyPi{\Cinfty\kern-3.5mm_{_{\bold\Pi}}\kern1.45mm}
\def\CinftyS{\Cinfty\kern-3.9mm_{_{\Cal S}}\kern1.45mm}
\def\wave{\hbox{\font\†=cmsy10\†\hbox{\char'164}\kern-2.35mm\hbox{\char'165}\kern.4mm}}
\def\wavee{\hbox{\font\†=cmsy8\†\hbox{\char'164}\kern-2.0mm\hbox{\char'165}\kern.4mm}} 
\def\barmj{\kern.25mm\bar{\hbox{\font\=cmr10\\char'021}}\kern.4mm}

\def\sigrd{\sigma\kern-.3mm_{_{rd}}\kern.15mm} 
\def\taurd{\tau_{_{rd}}\kern.15mm}
\def\tsigrd{\tau\sigma\kern-.3mm_{_{rd}}\kern.15mm} 
\def\tauR#1{\tau_{_{I\!\!R}}\kern-1.5mm^{#1}}
\def\tauC{\tau\kern-.3mm\lower.9mm\hbox{\font\=cmmi8\c}\kern.35mm}
\def\RN{I\!\!R\kern.3mm^{\hbox{\font\=cmmi6\N}}} 
\def\QTN{Q\kern.1mm_{\lower.2mm\hbox{\font\=cmmi6\T}}^{\kern.2mm\hbox{\font\=cmmi6\N}}} 

\def\leLCS-{{\le}{}_{_{\roman{LCS}}}\text{\sp-\sp}}

\def\sixroman#1{\hbox{\font\=cmr6\#1\kern.1mm}}
\def\sNor#1{\kern.25mm\lower.38mm\hbox{$_{#1}$}}
\def\subtext#1{\raise.2mm\hbox{$_{_{\kern0.15mm\text{#1}}}$}}
\def\Centerline#1\par#2\par#3{\noindent#1\phantom{#3}\hfill#2\hfill\phantom{#1}#3}

\def\usimW{\kern.35mm\lower.2mm\hbox{$_{_\sim}$}\kern-2.6mmW\raise.1mm\hbox{$^{^{}}$}}
\def\ubarW{\kern.4mm\underline{\kern-.4mmW\kern-1.7mm}\kern2.2mm}
\def\uvarPi{\kern.15mm\underline{\kern-.15mm\varPi\kern-.85mm}\kern.85mm}%
\def\uOmega{\kern.3mm\underline{\kern-.3mm\Omega\kern-.3mm}\kern.3mm}%
\def\CW{C\kern-.35mm\lower.1mm\hbox{$_{_W}$}\kern-2.1mm{^{}}}
\def\Topma{\roman{{Top_{}}_{ma}\kern.15mm}}

\def\overbar#1{\kern.65mm\overline{\kern-.6mm#1\kern-.6mm}\kern.65mm} 

\def\prodc{\prod{_{_{\kern-.3mm\bold c\kern.15mm}}}}

\def\prodmea{\kern1mm\raise.5mm\hbox{\font\=cmbsy5\\char'012}\kern1mm}
\def\sigmaA{\sigma_{\kern-0.4mm_A\kern.2mm}}
\def\loint{\int\kern-1.2mm\lower.45mm\hbox{$\underline{\phantom.}$}\kern1mm}
\def\upint{\kern.5mm\raise1.85mm\hbox{$\overline{\phantom.}$}\kern-1.2mm\int}
\def\Pint{\lower.7mm\hbox{$^{^{\roman P}}$}\kern-1.4mm\int}
\def\weakint{\lower.3mm\hbox{$^{^{\roman w}}$}\kern-1.7mm\int}
\def\sLL#1#2{\lower.8mm\hbox{{\font\=cmmi6\L}\raise1mm\hbox{\font\=cmmi5\#1}\kern.2mm{\font\=cmr6\\char'050}{\font\=cmmi6\#2}{\font\=cmr6\\char'051}}}


\hsize125mm \vsize253truemm \parskip.5mm \overfullrule=0mm 
\hoffset30truemm\voffset8truemm 
\hcorrection{-1truein}\vcorrection{-1.15truein}


\document \baselineskip4.5mm \newcount\firstpage\firstpage=1\pageno=\firstpage
\RunMyHead{S.\ Hiltunen}{}{}{Banach spaces over bounded domains}

{\font\=cmss12 scaled\magstep1\\centerline{%
                On the definition of some Banach spaces over}\vskip1.5mm\centerline{%
                  bounded domains with irregular boundary}}\vskip3mm\centerline{\font\=cmr12\%
                               by}\vskip2mm\centerline{\font\=cmmi12\%
                      Seppo\kern1mm Hiltunen}\vskip5mm\noin
\abstract{{\bf Abstract. }This note aims to clarify the interrelations of
  certain inequivalently defined Banach spaces denoted by $
C^{\ssp i}(\overline\Omega)$ for a natural number $i$ and a bounded open set $
\Omega\ssp$. We give some sufficient conditions for the equality of these
spaces, and present examples to show that the spaces indeed can be unequal for
$\Omega$ having irregular boundary.\vskip1mm\noin{\bf%
Keywords: }Banach space, differentiable function, irregular boundary.\vskip1mm\noin{\bf%
Subject classification:} 46E15. (AMS 2000)} 


\subhead

                          The setting and results                         

For fixed $i\ssp,\smb N\in\N=\{\ssp 1\ssp,2\ssp,\ldots\,\}$ and a bounded open
$\Omega$ with closure $Q$ in $\biit R\,\yi N$, we consider three Banach\sp
((iz)\sp able) locally convex spaces for which we generally have $
C^{\ssp i}(\overline\Omega)\overset\iota\to\hookleftarrow C^{\ssp i}(Q)\le
C^{\ssp i}(\overline\Omega)\subtext{\sp H\"o}\ssp$. Here $\iota$ defined by $x
\mapsto x\,|\,\Omega$ is a {\it strict morphism\ssp} in the sense of
[\,Ho\ssp; Def.\ 2.5.1, p.\ 100\,]\ssp, which means that $\iota$ is a linear
homeomorphism onto the subspace $\rng\iota$ equipped with the induced
topological vector space structure. Our concern is under what additional
assumptions $\iota$ is surjective or we also have $C^{\ssp i}(Q)=
C^{\ssp i}(\overline\Omega)\subtext{\sp H\"o}\ssp$.

We shall follow the conventions of [\,Hi\,]\ssp, from which we in particular
recall that topological (locally convex) vector spaces are understood to be
pairs $E=(X\sp,\cal T\sp)\ssp$, where $X=\sigrd E$ is the underlying vector
(space )\ssp structure, and $\cal T=\taurd E$ is the topology. The {\it
underlying set\ssp} of $E$ is $\vecs E=\bigcup\sp\cal T_{\sp}$, and the zero
vector is $\bnull E\ssp$. Hence contrary to the usual customs, we shall
systematically use distinct notations for the underlying set and the
topological vector space structure. For example, the set of real numbers is $
\Re=\vecs\biit R\,$. For locally convex $E$ and $F\sp$, the notation $E\le F$
expresses the fact that the identity $\idv F$ is a continuous linear map $F\to E\ssp$.

For the precise formal definition of $F=C^{\ssp i}(Q)$ referring the reader to
see [\,Hi\ssp; 2.3 Diff., p.\ 9\,]\ssp, roughly the definition is as follows.
The underlying set $\vecs F$ has as members exactly the functions $x:Q\to\Re$
which on the interior $\roman{Int_{\sp}}_{{\tau_{}}_{rd}\sp\bmii{7}R^N\sp}Q$
posses iterated partials $\partial^{\,\alpha}x$ for $|\ssp\alpha\sp|\le i$
having continuous extensions $x^{\,\alpha}$ on $Q\ssp$. The topology $\taurd F
$ is that of uniform convergence of all these $x^{\,\alpha}$ on $Q\ssp$.
Observe that we have $\Omega\inc
\roman{Int_{\sp}}_{{\tau_{}}_{rd}\sp\bmii{7}R^N\sp}Q$ and that $\Omega\not=
\roman{Int_{\sp}}_{{\tau_{}}_{rd}\sp\bmii{7}R^N\sp}Q$ is possible for example
if $\Omega=Q\sn\setminus\sn S$ for a closed $S$ with no interior points.

In the books [\,A\,]\ssp, [\,T\,]\ssp, and [\,W\,]\ssp, complex valued
functions are considered but for simplicity, we shall here modify the
definitions so as to produce spaces of real functions. In [\,T\ssp; p.\
xv\,]\ssp, the Banach space $E=C^{\ssp i}(\overline\Omega)$ is defined so that $
\vecs E$ has as members exactly the functions $x:\Omega\to\Re$ having partials
$\partial^{\,\alpha}x$ for $|\ssp\alpha\sp|\le i$ on $\Omega$ which can be
extended to continuous $x^{\,\alpha}$ on $Q\ssp$. The topology $\taurd E$ is
that of uniform convergence of all $\partial^{\,\alpha}x$ on $\Omega\ssp$. One
easily sees that this is equivalent to the uniform convergence of $
x^{\,\alpha}$ on $Q\ssp$. This proves the assertion that $\iota:F\to E$ is a
strict morphism. One easily verifies the implication $\,\Omega=
\roman{Int_{\sp}}_{{\tau_{}}_{rd}\sp\bmii{7}R^N\sp}Q\imply\rng\iota=\vecs E
\ssp$. That $\rng\iota\not=\vecs E$ is possible is shown in

\NSN{\bf
1} \Examplee We construct a $\tauR 2\,$--\,connected open $\Omega$ and $x\in
  \bigcap\ssp\{\,\vecs C^{\ssp i}(\overline\Omega):i\in\N\,\}$\Newline such
that no $\bar x\in\vecs C^{\ssp 1}(Q)$ with $x\inc\bar x$ exists. Let $I=
[\,0\,,1\,]$ and $A=\{\,\sigma\ssp(\biit a):\biit a\in$\Newline $\{\ssp 0\,,
2\ssp\}^{\ssp I\!\!N}\ssp\}\ssp$, where $\sigma\ssp(\biit a)=
\sum_{\ssp n\sp=\sp 1}^{\,\infty}a{_{}}_n\sp 3^{\sp-n}$ for $\biit a=\seq{\,
a{_{}}_n\ssp}\ssp$. Then $A$ is the $\tau_{_{I\!\!R}}\ssp$--\,compact\Newline
"Cantor set\ssp" of measure zero, for which in [\,R\ssp; 7.16\ssp(b)\ssp, p.\
145\,] 
       is constructed a continuous \hfill surjection \hfill $\varphi:I\to I$ \hfill
which \hfill is \hfill differentiable \hfill at \hfill every \hfill $s\in
I\setminus\sn A$ \hfill with\linebreak $\varphi\ssp'(s)=\varphi\ssp(0)=0\,$.
Taking $O=\{\ssp\seq{\ssp s\ssp,t\ssp}:-1<s\ssp,t<1\,\}\ssp$, we put $\Omega=
O\sn\setminus\sn S\ssp$, where $S=\{\ssp\seq{\ssp s\ssp,t\ssp}:s\in A$ and $t
\in I\,\}\ssp$. Let $x=\bar x\,|\,\Omega$ for $\bar x:Q\to I$ defined by $
x\ssp(\eta)=0$ for $\eta=\seq{\ssp s\ssp,t\ssp}\in Q$ with not $0<s\ssp,t\le 1
\ssp$, and $x\ssp(\eta)=\varphi\ssp(s)\,e^{\sp-t^{-1}}$ for $\eta\in Q$ with $
0<s\ssp,t\le 1\ssp$. It is an easy exercise for the reader to verify our
assertion above.                                                          \NS

In [\,W\ssp; p.\ 2\,]\ssp, the definition of a space $C^{\ssp i}(\bar\Omega)$
is expressed a bit vaguely by requiring the functions together with the
partials to "be continuous on $\bar\Omega\ssp$". We shall interpret this to
mean that they have continuous extensions, and then Wloka's definition equals
that of Treves, and hence $C^{\ssp i}(\bar\Omega)=E\ssp$. In [\,A\ssp; 1.26,
p.\ 9\,]\ssp, a space $C^{\ssp i}(\bar\Omega)\subtext{Ad}$ is defined so as to
have its underlying set formed by functions $x:\Omega\to\Re$ which together
with the appropriate partials are bounded and {\it uniformly continu- ous\ssp}
on $\Omega\ssp.$ The topology is that of uniform convergence of all these on $
\Omega\ssp$. Since we assume $\Omega$ to be bounded, uniform continuity is
equivalent to possessing a con- tinuous extension to $Q\ssp$, and hence we see
$C^{\ssp i}(\bar\Omega)\subtext{Ad}=E$ to hold.

We have now seen that the definitions of Adams, Treves, and Wloka all specify
the same space $E=C^{\ssp i}(\overline\Omega)\ssp$, which we consider the {\it
usual\ssp} one. By example 1\ssp, the set $\vecs E$ may have as members $
x\,|\,\Omega$ with $x:Q\to\Re$ rather pathological.

In a footnote in [\,H\"\ssp; p.\ 190\,]\ssp, the notation "$
C^i(\overline\Omega)\ssp$" is introduced to mean the {\it set\ssp} $
S\subtext{H\"o}$ of all functions $x:Q\to\Re$ possessing some $\bar x\in\vecs
C^{\ssp i}(\Re\,\yi{N\sp})$ with $x\inc\bar x\ssp$.
                                                    However, there is no
specification for a vector structure nor a topology to make this set the
underlying one of a Banach space. We now wish to construct a Banach space $H=
C^{\ssp i}(\overline\Omega)\subtext{\sp H\"o}$ with $\vecs H=S\subtext{H\"o}
\ssp$. A natural route is the following.
                                        
We let $H$ be the top\,-\,linear isomorphic image of $G\sp/\sp N$ under $x+N
\mapsto x\,|\,Q\ssp$, when $G$ and $N$ are as follows. We take $G=
C^{\ssp i}\KN{1.6}\LHB{.2}{\ai{bd}}\sp(\Re\,\yi{N\sp})\ssp$, the Banach space
of $C^{\ssp i}$ functions $\Re\,\yi N\to\Re$ having all partials up to order $
i$ bounded, and equipped with the topology of uniform convergence of all these
partials. We let $N$ be the linear subspace formed by all $x\in\vecs G$ having
$x\sp\image\sn Q=\{\sp 0\sp\}\ssp$. Since $\taurd G$ is stronger than the
topology of pointwise convergence, we see that $N$ is $\taurd G\,$--\,closed,
and hence $H$ is a Banach space since by [\,Ho\ssp; Thm.\ 2.9.2, p.\ 138\,] or
[\,Jr\ssp; Prop.\ 4.4.1, p.\ 80\,] the quotient of any complete metrizable
topological vector space by a closed linear subspace is complete.
                                                                  To see that
with this definition we indeed have $\vecs H=S\subtext{H\"o}\ssp$, we only
need to observe that for any $x\in\vecs C^{\ssp i}(\Re\,\yi{N\sp})\ssp$, we
have $y=\chi\,x\in\vecs G$ with $x\,|\,Q=y\,|\,Q$ for any compactly supported
smooth $\chi$ taking the value $1$ on $Q\ssp$.

Observe that the space $H=C^{\ssp i}(\overline\Omega)\subtext{\sp H\"o}$ only
depends on $Q=\roman{Cl_{\sp}}_{{\tau_{}}_{rd}\sp\bmii{7}R^N\sp}\Omega$ unlike
$C^{\ssp i}(\overline\Omega)$ whose specification requires the knowledge of $
\Omega\ssp$. We easily see $C^{\ssp i}(Q)\le H$ to hold. Neglecting the
problem whether generally even $H$ is a topological linear subspace of $
C^{\ssp i}(Q)\ssp$, we proceed to give in Proposition 2 below an additional
sufficient condition for $H=C^{\ssp i}(Q)$ to hold. For this, we agree to say 
that the set $Q$ {\it has a $C^{\ssp i}$ boundary\ssp} if{}f for every $\eta
\in{\partial_{\ssp}}_{{\tau_{}}_{rd}\sp\bmii{7}R^N\sp}Q$ there is a $
C^{\ssp i}$ diffeomorphism $\phi:\mathbreak B_{\bmii7R^N}(1)\to U$ with $U$ a
neighborhood of $\eta$ satisfying $\phi\,[\,
                                            B^{_{^{\,+}}}_{\bmii7R}\LHB{.2}{_{\sn^N}}(1)
                                                        \,]=U\cap Q\ssp$. \hfill
Here we have the open Euclidean ball $B_{\bmii7R^N}(1)=\{\,\eta\in\Re\,\yi N\sn
:|\sp\eta\sp|<1\,\}$ and and its subset $
                                         B^{_{^{\,+}}}_{\bmii7R}\LHB{.2}{_{\sn^N}}(1)
                                                     $ having as members
exactly the $\eta=\seq{\sp t\sp}\conc\bar\eta\in B_{\bmii7R^N}(1)$ with $t\ge
0\,$. Recall from [\,Hi\ssp; p.\ 5\,] that here $\eta=\seq{\,t\ssp,
\eta{_{}}_1\sp,\ldots\,\eta{_{}}_{\ssmb N\sp-\sp 1}\ssp}$ if $\bar\eta=\seq{\,
\eta{_{}}_1\sp,\ldots\,\eta{_{}}_{\ssmb N\sp-\sp 1}\ssp}\ssp$.

\newProCla 2 Proposition.

Let $i\ssp,\smb N\in\N\sp$, and let $\Omega$ be a bounded open set in $
\biit R\,\yi N$ with clos- ure $Q=
\roman{Cl_{\sp}}_{{\tau_{}}_{rd}\sp\bmii{7}R^N\sp}\Omega$ having a $
C^{\ssp i\,}$--\,boundary $C\ssp$. Then $
C^{\ssp i}(\overline\Omega)\subtext{\sp H\"o}=C^{\ssp i}(Q)\ssp$.

\Prooff Writing $H=C^{\ssp i}(\overline\Omega)\subtext{\sp H\"o}$ and $F=
C^{\ssp i}(Q)\ssp$, we have the Banach spaces $H$ and $F$ for which we already
saw $F\le H\sp$. If $\vecs F\inc\vecs H\sp$, the open mapping theorem shows $H
\le F\sp$, and hence $H=F\sp$. \hfil \hbox{\,Consequently, \,arbitrarily \,
fixing $\,x\in\vecs F\sp$, \,it}\linebreak suffices to show $x\in\vecs H\sp$,
i.e., to find a $C^{\ssp i}$ function $\bar x:\Re\,\yi N\to\Re$ with $x\inc\bar x\ssp$.

For the proof of existence of $\bar x\ssp$, we first show that given any $u\in
\vecs C^{\ssp i}(
                 B^{_{^{\,+}}}_{\bmii7R}\LHB{.2}{_{\sn^N}}(1)
                             )\ssp$, there is $\bar u\in\vecs C^{\ssp i}(
B_{\bmii7R^N}(1))$ with $u\inc\bar u\ssp$. Indeed, letting $\seq{\,
a{_{}}_0\ssp,\ldots\,a{_{}}_i\ssp}$ solve the linear system $
\sum_{\,l\ssp=\sp 1\sp}^{\,i\ssp+\sp 1}(-l\sp)^{\sp-j}a{_{}}_{l\ssp-\sp 1}=1$
for $j\in i\ssp\yplus$, \hfill for $\eta=\seq{\sp t\sp}\conc\bar\eta$ with $
t<0\,$, it suffices to define \hfill $\bar u\ssp(\eta)=
\sum_{\,l\ssp=\sp 1}^{\,i\ssp+\sp 1}a{_{}}_{l\ssp-\sp 1}\sp
u\ssp(\seq{-l^{\sp-1\sp}t\,}\conc\bar\eta\sp)\ssp$. \hfill See [\,W\ssp;
p.\ 101\,] for the computational\linebreak detail. Using this result together
with the assumption that $Q$ has $C^{\ssp i}$ boundary, we see that every $
\smb P\in C$ has some $C^{\ssp i}$ function $y$ defined on an open
neighborhood of $\smb P$ with $y\,|\,Q\inc x\ssp$. We may here take $y=\bar u
\circ(\sp\phi_{}^{\sp-\iota\sp})\ssp$, where $u=x\circ\phi\,$.

Finally, by compactness of $C\sp$, we have some finite sequence $
\seq{\,x{_{}}_0\ssp,\ldots\,x{_{}}_n\ssp}$ with $C\inc\bigcup\cal A$ of these
local extensions, and a smooth partition $\seq{\,\chi{_{}}_\nu:\nu\in N\,}$ of
unity subordinate to $\cal A=\{\,\roman{dom\,}x{_{}}_i:i=0\ssp,\ldots\,n\,\}
\ssp$. \hfill Writing $\,U=\Re\,\yi N\sn\setminus\sn(\sp Q\sn\setminus\sn
\bigcup\cal A\sp)\ssp$, \,and taking $i{_{}}_\nu=\min\ssp\{\,i\in n\sp\yplus\sn
:\roman{supp\,}\chi{_{}}_\nu\inc\roman{dom\,}x{_{}}_i\ssp\}$ for all $\nu\ssp
$, we define the $C^{\ssp i}$ functions $y{_{}}_\nu$ on $U$ by $\eta\mapsto
\chi{_{}}_\nu\sp(\eta)\,x_{i_\nu}(\eta)$ for $\eta\in\dom x_{i_\nu}\ssp$, and
$\eta\mapsto 0$ otherwise. \hfill Putting $y=\sum_{\ssp\nu\ssp\in\sp N\sp}
y{_{}}_\nu\ssp$, we have $y\,|\,Q\inc x\ssp$, whence with $\bar x=x\cup y$ we
are done.                                                             \newQED

By the discussion in [\,He\ssp; pp.\ 186\,--\,187\,]\ssp, it is obvious that
the sufficient condition for $C^{\ssp i}(\overline\Omega)\subtext{\sp H\"o}=
C^{\ssp i}(Q)$ given by Proposition 2 is not necessary. Omitting the task of
formulating a more general (and more complicated) sufficient condition, we
next proceed to give an example showing that indeed $
C^{\ssp i}(\overline\Omega)\subtext{\sp H\"o}\not=C^{\ssp i}(Q)$ is possible
if $Q$ has sufficiently irregular boundary.

A disconnected $\Omega\inc\Re$ for $i=1$ can easily be constructed from the
idea in [\,Hi\ssp; Remarks 2.4, p.\ 10\,]\ssp. Namely, we let $\Omega$ be the
interior of $Q=\bigcup\ssp\{\ssp I{_{\sn}}_n:n\in\No\ssp\}\ssp$, where $
I{_{}}_0=[-1\ssp,0\,]$ and $I{_{\sn}}_n=[\,s{_{}}_n\ssp,
\frac32\ssp s{_{}}_n\ssp]$ with $s{_{}}_n=2^{\sp-n}$ for $n\in\N\sp$. Putting
$x=$ $(\sp\roman{id\,}I{_{}}_0)\cup\bigcup\ssp\{\ssp\seq{\,s-s{_{}}_n\sn:s\in
I{_{\sn}}_n\ssp}:n\in\N\ssp\}\ssp$, \hfill one verifies that we have $x\in
\vecs C^{\ssp 1}(Q)$ and $x\not\in\vecs
C^{\ssp 1}(\overline\Omega)\subtext{\sp H\"o}\ssp$. \hfill From the same idea,
we now wish to establish a connected $\Omega$ in dimension two with the
analogous property in

\NSN{\bf
3} \Examplee There are $\Omega\in\tauR 2=\taurd\biit R^{\,2}$ and $x$ such
  that $\Omega$ is $\tauR 2\ssp$--\,connected, and with $Q=
\roman{Cl_{\sp}}_{\tau_{I\!\!R}^{\sp 2}\sp}\Omega\ssp$, we have $x\in\vecs
C^{\ssp 1}(Q)$ and $x\not\in\vecs
C^{\ssp 1}(\overline\Omega)\subtext{\sp H\"o}\ssp$. \hfill Indeed, we can\linebreak
take $\Omega=\roman{Int_{\sp}}_{\tau_{I\!\!R}^{\sp 2}\sp}Q\ssp$, where $Q$ and
$x$ are constructed as follows. With $b{_{}}_n=2^{-n}$ and $a{_{}}_n=
\frac34\ssp b{_{}}_n\ssp$, \hfill writing \hfill $A{_{}}_n=\{\ssp\seq{\ssp
s\ssp,t\ssp}:a{_{}}_n\le s\le b{_{}}_n$ and $0<t\le 1\,\}\ssp$, \hfill and \hfill
$B=\mathbreak\{\ssp\seq{\ssp s\ssp,t\ssp}:-1\le s\ssp,t\le 1$ and not $
0<s\ssp,t\,\}\ssp$, \hfill we \hfill obtain \hfill the \hfill "comblike\ssp" \hfill
set \hfill $Q=\mathbreak\bigcup\ssp\{\,A{_{}}_n\sn:n\in\No\ssp\}\cup B\ssp$. \hfill
We now take $x=\bigcup\ssp\{\,x{_{}}_n\sn:n\in\No\ssp\}\cup(\sp\chi\,|\,B\sp)
\ssp$, where with $\chi=\{\ssp(\seq{\ssp s\ssp,t\ssp}\ssp,s\,t^{\,2\sp}):
s\ssp,t\in\Re\,\}\,$, and also writing $c{_{}}_n=b{_{}}_n-a{_{}}_n\,$, \hfill
we define\,%
\footnote{note the definition $\{\,\frak t\ssp(\sp n\sp,s\ssp,t\sp):n:
  \frak P\sp(\sp n\sp,s\ssp,t\sp)\ssp\}=\{\,z:\exi{s\ssp,t}\,z=
  \frak t\ssp(\sp n\sp,s\ssp,t\sp)$ and $\frak P\sp(\sp n\sp,s\ssp,t\sp)\ssp\}
  $ obtained as a particular case of the general definition schema $\{\,
  \frak T:\frak x_{\sp l_2+1},\ldots\,\frak x_{\sp l_3}:\frak F\,\}=\{\,
  \frak x_{\sp 0}:\exi{\frak x_{\sp l_1}\sp,\ldots\,\frak x_{\sp l_2}}\,
  \frak x_{\sp 0}=\frak T$ and $\frak F\,\}\ssp$, when the term $\frak T$ has
  free variables $\frak x_{\sp 1}\sp,\ldots\,\frak x_{\sp l_3}$ and the
  formula $\frak F$ has free variables $\frak x_{\sp l_1}\sp,\ldots\,
  \frak x_{\sp l}\ssp$. Here we understand that $\frak x_{\sp i}$ and $
  \frak x_j$ are distinct variable\sp(symbol)\sp s unless $i=j\ssp$, and also
  that $1\le l_1\le l_2\le l_3\le l\ssp$. In the case where $l_2=l_3\ssp$, we
  abbreviate $\{\,\frak T:{\,}:\frak F\,\}=\{\,\frak T:\frak F\,\}\ssp$.
  Compare to the discussion in [\,Ky\ssp; 4 Notes, p.\ 6\,]\ssp.}
$x{_{}}_n=\mathbreak\{\ssp(\seq{\ssp s\ssp,t\ssp}\ssp,c{_{}}_n\ssp\chi\value
\seq{\,c_{\sp n}^{\sp-1}(\sp s-a{_{}}_n)\ssp,t\ssp}):n:\seq{\ssp s\ssp,t\ssp}
\in A{_{}}_n\ssp\}\ssp$.

Writing $\alpha\ar 1=\seq{\ssp 1\ssp,0\,}$ and $\alpha\ar 2=\seq{\,0\,,1\ssp}
\ssp$, to prove $x\in\vecs C^{\ssp 1}(Q)\ssp$, the only nontrivial task is to
verify continuity of $x$ and $
                              x^{\,\alpha_i}$ at points $\xi\ar 0=\seq{\,0\,,
t\ssp}$ for $0\le t\le 1$ and $\xi\ar 1=\seq{\ssp s\ssp,0\,}$ for $a{_{}}_n\le
s\le b{_{}}_n\ssp$. This we leave as a routine exercise for the reader. To\linebreak
prove \hfill $x\not\in\vecs C^{\ssp 1}(\overline\Omega)\subtext{\sp H\"o}\ssp
$, \hfill it \hfill suffices \hfill to \hfill observe \hfill that \hfill in \hfill
the \hfill opposite \hfill case \hfill we \hfill would\linebreak get $\,1=
\underset{s\to 0^-}\to\lim\partial\ar 1\sp x\value\sn\seq{\ssp s\ssp,1\ssp}
=\underset{n\to+\infty}\to\lim\,a_{\sp n}^{\sp-1\sp}(\sp x\value\sn\seq{\ssp
a{_{}}_n\sp,1\ssp}-x\value\sn\seq{\,0\,,1\ssp})=0\,$, \,a contradiction.

\NSN{\bf
4} \Remarkk The construction via which we obtained the space $
  C^{\ssp i}(\overline\Omega)\subtext{\sp H\"o}$ can be generalized as
follows. Assume that for $\varPi\sp,F\in\roman{LCS}\ssp(\bold R)$ with $\varPi
$ finite dimensional we are given a space $E=S\ssp(\varPi\sp,F\sp)\in
\roman{LCS}\ssp(\bold R)$ such that $\idv E\in\cal L\ssp(E\ssp,
F^{\,\upsilon_s\varPi\sp})\ssp$. For $Q\inc\vecs\varPi$ we then define \ $H=
S\ssp(\overbar Q_\varPi\sp,F\sp)=$\par\centerline{$
\leLCS-\sup\ssp\{\,G:F^{\,Q}\le G$ and $\seq{\,x\,|\,Q:x\in\vecs E\,}
                                 \in\cal L\ssp(E\ssp,G\sp)\ssp\}\,$.}

\noin Then $\seq{\,x\,|\,Q:x\in\vecs E\,}$ is a strict morphism $E\to H$ in
Horv\'ath's second sense [\,Ho\ssp; Def.\ 2.5.2, p.\ 106\,]\ssp, and hence a
linear homeomorphism $E\ssp/\sp N\to H$ is given by $x+N\mapsto x\,|\,Q$ when
we take $N=\{\,x\in\vecs E:x\image Q\inc\{\sp\bnull F\}\sp\}\ssp$. Observe
that the requirement $\idv E\in\cal L\ssp(E\ssp,F^{\,\upsilon_s\varPi\sp})$
guarantees that $N$ is $\taurd E\,$--\,closed, whence indeed we have $H\in
\roman{LCS}\ssp(\bold R)\ssp$. Since a quotient of any Banach space by a
closed subspace is Banach, we see that if $E$ is Banach, then so is $H\sp$.

The above general construction for example gives a meaning to the real Banach
space ${H^{\sp}}^s(Q)_{_{\bmii5C}}={H^{\sp}}^s(\sp\overbar Q_{\bmii6R^N},
\biit C_{_{I\!\!R}})$ of continuous complex valued "Fourier\,--\,wise defined
Sobolev\ssp" functions for any set $Q\inc I\!\!R\,\yi N$ when $\frac12\ssp
\smb N<s\in\Re\,$. 


\vskip5mm\centerline{\bf%
                             References}\vskip2mm
{\eightpoint\parskip.5mm\baselineskip3.5mm\leftskip9mm\parindent-9mm
   \def\item#1=#2..{[\,#1\,]\vskip-4mm\noindent#2.\vskip.5mm}
\item A=Adams, R.\ A.: {\it Sobolev Spaces\ssp}, Academic Press, New York -
  San Francisco - London 1975..
\item He=Hestenes, M.\ R.: Extension of the range of a differentiable
  function, {\it Duke Math.\ J\ssp}.\ Vol.\ 8 (1941) 183--192..
\item Hi=Hiltunen, S.: Differentiation, implicit functions, and applications
  to generalized well-posedness, {\it preprint\ssp} http://arXiv.org/abs/math.FA/0504268..
\item Ho=Horv\'ath, J.: {\it Topological Vector Spaces and Distributions\ssp},
  Addison--Wesley, Reading 1966..
\item H\"o=H\"ormander, L.: {\it Linear Partial Differential Operators\ssp},
  GMW {\bf116}, Springer, Berlin - G\"ottingen - Heidelberg 1963..
\item Jr=Jarchow, H.: {\it Locally Convex Spaces\ssp}, Teubner, Stuttgart 1981..
\item Ky=Kelley, J.\ L.: {\it General Topology\ssp}, GTM {\bf27}, Springer, 
  New York 1985..
\item R=Rudin, W.: {\it Real and Complex Analysis\ssp} ($3^{\roman{rd}}$
  ed.)\ssp, McGraw-Hill, New York 1987..
\item T=Treves, F.: {\it Basic Linear Partial Differential Equations\ssp},
  Academic Press, New York - San Francisco - London 1975..
\item W=Wloka, J.: {\it Partial Differential Equations\ssp}, Cambridge Univ.\
  Press, Cambridge 1992..
\par}\vskip3mm

{\eightpoint\baselineskip3.5mm\parskip0mm
  Seppo Hiltunen\vskip.5mm

  Helsinki University of Technology\par
  Institute of Mathematics, U311\par
  P.O.\ Box 1100\par
  FIN-02015 HUT\vskip.5mm
  FINLAND\vskip0mm
  e-mail: shiltune\,\@\,cc.hut.fi}

\enddocument